\newif\ifsmfart
\IfFileExists{smfart.cls}
   {\documentclass[12pt,english]{smfart}
    \IfFileExists{smfenum.sty}{\usepackage{smfenum}}{}
    \usepackage{bull}
    \smfarttrue}
{\message{^^J*** It would be better to typeset this file with smfart.cls ***^^J^^J}

\documentclass[12pt]{amsart}}
  
\usepackage{amssymb}
\usepackage{epsfig}
\usepackage{graphicx}
\numberwithin{equation}{section}

\input xy
\xyoption{all}

\advance\headheight 2pt
\calclayout
\allowdisplaybreaks[3]

\theoremstyle{plain}
\newtheorem{prop}[subsection]{Proposition}

\newtheorem{thm}[subsection]{Theorem}
\newtheorem{coro}[subsection]{Corollary}

\newtheorem{lemm}[subsection]{Lemma}

\theoremstyle{definition}

\newtheorem{ques}[subsection]{Question}
\newtheorem{conj}[subsection]{Conjecture}

\theoremstyle{remark}
\newtheorem{rem}[subsection]{Remark}

\newtheorem{exam}[subsection]{Example}

\newtheorem{nota}[subsection]{Notations}

\newcommand{\la}{\lambda}

\newcommand{\Q}{\Bbb Q}

\def\Bran{{\rm Bran}}
\def\Ram{{\rm Ram}}
\def\no{\noindent}

\def\cX{{\mathcal X}}

\def\cC{{\mathcal C}}

\def\cE{{\mathcal E}}

\def\cH{{\mathcal H}}
\def\cM{{\mathcal M}}

\def\cU{{\mathcal U}}

\def\ovl{\overline}

\def\ra{\rightarrow}

\def\F{{\mathbf F}}
\def\P{{\mathbb P}}
\def\Q{{\mathbb Q}}

\def\N{{\mathbb N}}
\def\cN{{\mathcal N}}

\def\sT{{\mathsf T}}

\def\Ram{{\rm Ram}}

\makeatother
\makeatletter

\author{Fedor Bogomolov}
\address{Courant Institute of Mathematical Sciences, N.Y.U. \\
 251 Mercer str. \\
 New York, NY 10012, U.S.A.}
\email{bogomolo@cims.nyu.edu}

\author{Yuri Tschinkel}
\address{Department of Mathematics \\
         Princeton University\\ 
         Fine Hall, Washington Road\\
         Princeton, NJ 08544-1000,  U.S.A.}
\email{ytschink@math.princeton.edu}

\title[Unramified correspondences]
{Unramified correspondences}

\begin{document} 
 
\date{\today}



\begin{abstract}
We study correspondences between algebraic curves 
defined over the separable closure of $\Q$ or $\F_p$.  
\end{abstract}

\maketitle

\tableofcontents

\setcounter{section}{1}
\section*{Introduction}
\label{sect:introduction}

A  class $\cC(\ovl{\Q})$ of complete 
algebraic curves over $\ovl{\Q}$ will be 
called {\em dominating} 
if for every algebraic curve $C'$ over $\ovl{\Q}$ 
there exist a curve $\tilde{C}\in \cC(\ovl{\Q})$ and  
a birational surjective
map $\tilde{C}\ra C'$. 
A curve $C$ will be called {\em universal} if the class 
$\cU_C(\ovl{\Q})$ of its unramified covers is dominating.

\begin{thm}[Belyi]
Every algebraic curve $C$ 
defined over a number field admits a surjective map onto
$\P^1$ which is unramified outside $(0,1,\infty)$. 
\end{thm}

In 1978 Manin pointed out
that Belyi's theorem implies 
the following

\begin{prop}
\label{prop:domi}
The class $\cM\cU(\ovl{\Q})$ consisting 
of modular curves and their unramified 
covers is dominating.
\end{prop}

There are many other classes of curves
with the same property, for example:
\begin{enumerate}
\item 
hyperelliptic curves and their unramified coverings;
\item 
the class $\cC\cU(\ovl{\Q}):=\cup_{n\in \N} \cC_n(\ovl{\Q})$, 
with $\cC_n(\ovl{\Q})$ consisting of curves  
with function field $\ovl{\Q}(z,\sqrt[n]{z(1-z)})$
and their unramified coverings. 
\item the class $\cC\cN(\ovl{\Q}):=\cup_{n\in \N} \cC\cN_n(\ovl{\Q})$ 
where $\cC\cN_n(\ovl{\Q})$ consists of all
unramified covers of any curve $C_n$ with the property that 
$C_n\ra\P^1$ is ramified in $(0,1,\infty)$ only and 
all local ramification indices of $C_n$ 
over 0 are divisible by 3, over 1 divisible by 2 and 
over $\infty$ divisible by $n$. 
In particular, we could take $C_n$ to 
be the modular curve $X(n)$.   
\end{enumerate}

\begin{proof}(Sketch) 
Let us consider the class of hyperelliptic curves and their
unramified covers. Let $C'$ be an arbitrary curve and 
$\sigma\,:\, C'\ra \P^1$ a {\em generic} map, branched over the points
$q_1,...,q_n$ (generic means that there is only one ramification point over
each branch point and all local ramification indices are equal to 2). 
Denote by $C$ a hyperelliptic curve whose 
ramification contains $q_1,...,q_n$.
Then $\tilde{C}:=C\times_{\P^1} C'$ 
is an unramified cover of $C$ which surjects onto $C'$.
For the classes $\cC\cU(\ovl{\Q})$ and  
$\cC\cN(\ovl{\Q})$ we use Belyi's theorem.
\end{proof}

\begin{ques}
\label{q:1}
Does there exist a universal algebraic curve $C$ (over $\ovl{\Q}$)?
\end{ques}

\begin{ques}
\label{q:2}
Does there exist a number $n\in \N$ such that
every curve defined over $\ovl{\Q}$ admits a surjective
map onto $\P^1$ with ramification over $(0,1,\infty)$
such that all local ramification indices are $\le n$? 
\end{ques}

\begin{ques}
\label{q:3}
Is every curve $C$ (over $\ovl{\Q}$) 
of genus $g(C)\ge 2$ universal?
\end{ques}

\begin{rem}
It is clear that an affirmative answer to Question~\ref{q:2} implies
a (constructive) affirmative answer to Question~\ref{q:1}.
\end{rem}

In this note we answer these questions in
a simple model situation: instead of $\ovl{\Q}$ 
we consider the (separable) closure $\ovl{F}_p$ of the finite field
${\F}_p$. 

\begin{thm}
\label{thm:main}
Let $p\ge 5$ be a prime and 
$C$ a hyperelliptic curve over 
$\ovl{{\F}}_p$ of genus $g(C)\ge 2$.
Then $C$ is universal: 
for any projective curve $C'$ there exist a finite 
\'etale cover $\tilde{C}\ra C$ and a  surjective regular map
$\tau\,:\,\tilde{C}\ra C'$. 
\end{thm}

\

In Section~\ref{sect:var} we prove the following geometric fact
(over arbitrary algebraically 
closed fields of characteristic $\neq 2,3$):

\begin{prop}
\label{prop:hypo}
Every hyperelliptic curve $C$ has a finite \'etale cover $\tilde{C}$ 
which surjects onto the genus 2 curve $C_0$ given by 
$\sqrt[6]{z(1-z)}$. 
In particular, if $C_0$ is universal 
then every hyperelliptic curve of genus $\ge 2$ is universal.
\end{prop}

\begin{rem}
Applying the Chevalley-Weil theorem we conclude   
that the Mordell conjecture (Faltings' theorem)
for $C_0$ implies the Mordell conjecture  
for every hyperelliptic curve of genus $\ge 2$.

The fact that there is some interaction between the arithmetic
of different curves has been noted previously. 
Moret-Bailly and Szpiro showed (see \cite{S}, \cite{M}) 
that the proof of an {\em effective} Mordell conjecture
for {\em one} (hyperbolic) 
curve (for example, $C_0$) implies the ABC-conjecture, 
which in turn implies an effective Mordell conjecture for {\em all}
(hyperbolic) curves
(Elkies \cite{E}). Here {\em effective} means an explicit bound
on the height of a $K$-rational point on the curve for all 
number fields $K$. Here again, Belyi's theorem is used in an essential way. 
\end{rem}

\noindent
{\bf Acknowledgments.} 
We have benefited from conversations with B. Hassett and 
A. Chambert-Loir.
The first author was partially 
supported by the NSF. The second author was partially
supported by the NSF and the Clay foundation.  

\

\section{Main construction}
\label{sect:2}

\begin{nota}
\label{nota:loci}
Let $\tau\,:\, C\ra C' $ be a surjective map of  algebraic curves. 
We denote by $\Ram(\tau)\subset C$ the ramification locus of $\tau$ and 
by $\Bran(\tau)=\tau(\Ram(C))\subset C'$ the branch locus of $\tau$.
For a point $q\in C$ we denote by $e_{q}(\tau)$ the local 
ramification index at $q$. 
We denote by 
$$
e(\tau):=\max_{q\in C}e_{q}(\tau)
$$ 
the maximum local ramification index of $\tau$. 
We say that $\tau$ has {\em simple} ramification if $e(\tau)\le 2$
and that $\tau$ is {\em generic} if in addition there is only one ramification 
point over each branch point. 
\end{nota}

\begin{rem}
Every curve admits a 
generic map onto $\P^1$, at least after a separable
extension of the ground field.
\end{rem}

Let $p\ge 5$ be a prime number. In this section 
we work over a separable closure $\ovl{\F}_p$
of the finite field $\F_p$.
First we show that there exists at least one curve satisfying 
the conclusion of Theorem~\ref{thm:main}.

\
 
\no
Let $\pi_0\,:\, E_0\ra \P^1$ be the elliptic curve given by
$$
\sqrt[3]{z(z-1)}.
$$
Let $\sigma_0\,:\, C_0\ra \P^1$ be the genus 2 curve given by
$$
\sqrt[6]{z(z-1)},
$$
and $\iota_0\,:\, C_0\ra E_0$ the corresponding 2-cover.
Clearly, $\iota_0$ has simple ramifications
over the preimages of $0,1$. 
Let $C$ be an arbitrary curve. Choosing a generic function
on $C$ we get a generic covering $\sigma\,:\, C\ra \P^1$ (such covering
is defined over $\ovl{\F}_p$). 
Assume further that $\Bran(\sigma)\subset \P^1$ 
does not contain $(0,1,\infty)$. 

\

Consider the diagram

\centerline{
\xymatrix{
  C \ar[d]_{\sigma} & \ar[l]\ar[d] C_1              & \ar[l]\ar[dd]C_2\\
\P^1       & \ar[l] \ar[d]_{\varphi} E_0   &                 \\
           &  E_0                          & \ar[l] C_0   
}
}

Here $C_1 = C \times_{\P^1} E_0$ 
(it is irreducible since $E_0\ra \P^1$ is a 2-cover). 
Then $C_1\ra E_0$ has simple 
ramification over a finite number  of points in $E_0$.
Recall that $E_0$ has a group scheme structure, and {\em all}
$\ovl{{\F}}_p$-points of $E_0$ are torsion points. 
This implies that there exists an \'etale map $E_0\ra E_0$
such that all ramification points of $C_1$ over $E_0$ are mapped to 
$0$. More precisely, any finite set of  
$\ovl{{\F}}_p$-points of $E_0$ 
is contained in the group subscheme $E_0^{et}[n]\subset E_0$ - 
the maximal \'etale subgroup 
of the multiplication by $n$-kernel $E_0[n]$
(for some $n\in \N$). For every positive integer $n$ there 
exists a positive multiple of $m$ of $n$
and an \'etale map
$E_0\ra E_0$ with kernel $E_0^{et}[m]$.

\

Taking the composition of $C_1\ra E_0$ with the multiplication 
by a suitable $m$, we get a (possibly new) surjective 
regular map $C_1\ra E_0$ which is ramified only 
over the zero point in $E_0$ and has the property  
that all the local ramification indices are at most 2.
Using this map let us define $C_2:=C_0\times_{E_0} C_1$. 
Consequently, any component of  $C_2$
surjects onto $C_1$ and is an \'etale covering of $C_0$
(ramification cancels ramification). 
This component satisfies the conclusion of Theorem~\ref{thm:main}.

\

\begin{lemm}
\label{lemm:ee}
Let $C$ be any smooth complete 
algebraic curve and $E$ any curve of genus 1. 
There exists a curve $C_1$ which surjects onto 
$C$ and $E$ such that the ramification of the map $C_1\ra E$ 
lies entirely over a single point of $E$ and
its local ramification indices are all equal to $2$. 
\end{lemm}

\begin{proof}
Consider a generic map $\sigma\,:\, C\ra \P^1$ 
with $e(\sigma)\le 2$. Choose a double cover $\pi\,:\, E\ra \P^1$
such that the branch loci 
$\Bran(\sigma)$ and $\Bran(\pi)$ on $\P^1$ are disjoint.  
Then the product $C_1:=C\times_{\P^1} E$ is
an irreducible curve which is a double cover of $C$.
The curve admits a surjective map $\iota_1\,:\, C_1\ra E$ 
with $e({\iota_1})\le 2$.
Similarly to the previous construction we can find an unramified
cover $\varphi\,:\, E\ra E$ such that the composition 
$\varphi\circ \iota_1\,:\, C_1\ra E$ is ramified only  over one point  
in $E$ and the local ramification indices are still equal to $2$. 
\end{proof}

\begin{coro}
\label{coro:ce}
Assume that some unramified covering $\tilde{C}$ of $C$ surjects
onto an elliptic curve $E$. Assume further that
there exists a point $q$ on $E$ such that {\em all} 
local ramification indices of the map $\tilde{C}\ra E$ 
over $q$ are divisible by 2. Then $C$ is universal.  
\end{coro}

\begin{proof}
It is sufficient to take the product of 
$\tilde{C}\times_E C_1$. Any irreducible component of the resulting 
curve will be an unramified covering of $\tilde{C}$ (and hence $C$)
and will admit a surjective map onto $C_1$ and $C$. 
\end{proof}

\begin{coro}[Theorem~\ref{thm:main}]
\label{coro:c}
Every hyperelliptic curve $C$    
over $\ovl{\F}_p$ (with $p\ge 5$)
of genus $\ge 2$ is universal.
\end{coro}

\begin{proof}
Consider the standard 
projection $\sigma\,:\, C\ra \P^1$ (of degree 2). 
Its branch locus ${\rm Bran}(\sigma)$ consists of $2g+2$ points. 
Let $\pi\,:\, E\ra \P^1$ be a double cover such that
${\Bran}(\pi)$ is contained in ${\Bran}(\sigma)$. Then the product
$\tilde{C}=C\times_{\P^1}E$ is an unramified double cover of $C$. 
Moreover, $\tilde{C}$ is a double cover of $E$ with ramification at most
over the preimages in $E$ of the points in 
${\Bran}(\sigma)\setminus {\Bran}(\pi)$. 
We now apply Corollary~\ref{coro:ce}. 
\end{proof}

In {\em finite} characteristic, there are many other  (classes of) 
universal curves. For example, cyclic coverings with ramification in 
3 points, hyperbolic modular curves, etc.  
Thus it seems plausible to formulate the following

\begin{conj}
\label{conj:main}
Any smooth complete curve $C$ of genus $g(C)\ge 2$ defined over 
$\ovl{{\F}}_p$ (for $p\ge 2$) is universal. 
\end{conj}

\section{The case of characteristic 0}
\label{sect:case}

In this section we work over $\ovl{\Q}$. We show that 
the method outlined in Section~\ref{sect:2} can employed
in characteristic zero to produce natural infinite sets of 
algebraic points on $\P^1$ which occur as
ramification points of surjective 
maps from $\P^1_2$ to $\P^1_1$ branched over $(0,1,\infty)\in \P^1_1$ only
and having an {\em a priori} bound on the ramification index
(here $\P^1_1$ and $\P^1_2$ are two different copies of the projective
line $\P^1$). 

\

Notice that, in principle, it is easy to produce {\em some} 
sets of points (of any finite cardinality)
with this property: Take an $n\ge 6$ and any triangulation of $\P^1_2$
with vertices of index $\le n$. A barycentric subdivision of each 
such triangulation defines a function from $\P^1_2$ to $\P^1_1$ with local
ramification indices $\le 2n$ (for more details see \cite{bh}). 
Therefore, any curve with bounded ramification over
this set of vertices will have bounded ramification over $\P^1_1$.
However, we have no explicit control over the coordinates of the ramification
points on $\P^1_2$.

\

An (obvious) analogous  
way to control ramification 
indices is to consider the following diagram

\

\centerline{
\xymatrix{
E\ar[d]_{\phi_n}\ar[r]^{\pi} & \P^1_2\ar[d]^{\varphi_{n,E}} \\
E            \ar[r]_{\pi} & \P^1_1, 
}}

\

\noindent
where the map $\phi_n$ is the quotient by the subscheme of $n$-torsion points
and the maps $E\ra \P^1$ are the standard double covers, ramified over 
$(0,1,\infty, \lambda)$. 
Clearly, all the ramification points of $\varphi_{n,E}$ (in $\P^1_2$)
are over $0,1,\infty$ and $\la$  (in $\P^1_1$) 
and $e(\varphi_{n,E})=2$. 
Belyi's theorem gives a map 
$\beta\,:\, \P^1_1\ra \P^1_0$, 
which ramifies only over the points
$(0,1,\infty)\in \P^1_0$, maps $\{0,1,\infty,\la\}\subset \P^1_1$ into 
$\{ 0,1,\infty\}\subset \P^1_0$ and has local ramification indices $\le n$.
Moreover, it provides an explicit bound on $\deg(\beta)$ and, consequently, 
on $e(\beta)$ (in terms of the absolute height of $\la$).  
Let $\beta_{\la}\,:\, \P^1_1\ra \P^1_0$ be a map such that
$$
e(\beta_{\la})=\inf_{\beta} \{  e_{\beta}\}
$$
over the set of all maps as above. 
Then the map $\beta_{\la}\circ \varphi_{n,E}\,:\, \P^1_2\ra \P^1_0$ 
ramifies over three points only 
and has index $e(\beta_{\la}\circ \varphi_{n,E}) \le 2n$.  
Let 
$$
R_E:=\pi(E(\ovl{\Q})_{\rm tors}) \subset \P^1_2(\ovl{\Q})
$$ 
be the image of the torsion points of $E$. 
Let $\sigma\, :\, C\ra \P^1_2$ be
any map ramified only in a subset of $R_E$.  
Let $\pi:=\beta_{\la}\circ \varphi_{n,E}\circ \sigma$.
Then 
$$
e(\pi)\le 2e(\sigma)\cdot e(\beta_{\la}).
$$

\

A natural application of the construction in Section~\ref{sect:2} is 
as follows:

\begin{exam}
\label{exam:3}
Let $\pi\,:\, E\ra \P^1$ be a triple cover 
with $\Bran(\pi)=\{ 0,1,\infty\}$
($E$ is a CM elliptic curve with $j$-invariant $0$). 
Consider the following diagram

\

\centerline{
\xymatrix{
      & E\ar[d]_{\phi_n}\ar[r]^{\pi} & \P^1_2\ar[d]^{\varphi_{n,E}} \\
C_0 \ar[r] & E\ar[r]_{\pi} & \P^1_1, 
   }}

\

\noindent
where $C_0$ is a curve of genus $g(C_0)=2$ given by $\sqrt[6]{z(z-1)}$, 
$\phi_n$ is the quotient map by the subscheme of torsion points of
order $n$, and $\varphi_{n,E}$ the corresponding map from $\P^1_2$ to $\P^1_1$
ramified only over $(0,1,\infty)$.
Let $\cX_g=\{ X\}$ be the subset of curves 
of genus $g$ admitting a map $\sigma_X\,:\, X\ra \P^1_2$ such that
\begin{itemize}
\item $e(\sigma_X)=2$; 
\item $\Bran(\sigma_X)\subseteq \pi(E(\ovl{\Q})_{\rm tors})$.
\end{itemize}
Then, for any $X\in \cX_g$ the map 
$$
\varphi_{n,E}\circ \sigma_X\,:\, X\ra \P^1_1
$$
has index $e(\varphi_{n,E}\circ \sigma_X)\le 6$ 
and there exists an unramified cover $\tilde{C}\ra C_0$
surjecting onto $X$.
Moreover, $\cX_g$ is {\em dense} (in real and $p$-adic topologies) 
in the natural Hurwitz scheme $\cH_g$
parametrizing curves of genus $g$. 
\end{exam}

The set of curves dominated by unramified covers of $C_0$ is much 
larger than $\cX_g$. 
Indeed, consider any 4-tuple of points 
in 
$$
\pi(E(\ovl{\Q})_{\rm tors})\subseteq \P^1_2
$$ 
and an elliptic curve $E'$ obtained
as a double cover of $\P^1_2$ ramified in those 4 points. 
Then $E'$ is also dominated by unramified covers of $C_0$ and we can iterate
the above construction for $E'$.

\section{Geometric constructions}
\label{sect:var}

Let $(E,q_0)$ be an elliptic curve, 
$q_1$ a torsion point of order two 
on $E$ and $\pi\,:\, E\ra \P^1$
the quotient with respect to the involution induced by $q_1$.    
Let $n$ be an odd positive integer and $\varphi_{n,E}\,:\, \P^1_2\ra \P^1_1$
the map induced by 

\

\centerline{
\xymatrix{
E\ar[d]_{\phi_n}\ar[r]^{\pi} & \P^1_2\ar[d]^{\varphi_{n,E}} \\
E            \ar[r]_{\pi} & \P^1_1.  
}}

\

\noindent
Any quadruple $r=\{r_1,...,r_4\}$ of four
distinct points in  $\varphi^{-1}_{n,E}(\pi(q_0))$
defines a genus 1 curve $E_r$ (the double cover of $\P^1$ ramified 
in these four points).

\begin{prop}
\label{prop:tree}
Let $\iota\,:\, C\ra E$ be any finite cover 
such that all local ramification indices
over $q_0$ are even. 
Then there exists an unramified cover 
$\tau_r\,:\, C_r\ra C$ which 
dominates $E_r$ and has only even local ramification indices
over some point in $E_r$. 
\end{prop}

\begin{proof}
Assume that  $n\ge 3$ and  consider 
the following diagram
\

\centerline{
\xymatrix{
C\ar[d]_{\iota} & \ar[l]_{\tau_2}\ar[d]^{\iota_2}C_2 
& \ar[l]_{\tau_r}\ar[d]^{\iota_r} C_r \\
E\ar[d]_{\pi}   & \ar[l]_{\varphi_n}\ar[d]^{\pi} E   
&                     \ar[d]^{\pi_r} E_r   \\
\P^1_1            & \ar[l]_{\phi_{n,E}} \P^1_2  & \P^1_2,                    
   }}

\

\noindent
where $E_r$ is a double cover of $\P^1_2$ ramified in any quadruple of points
in the preimage $\phi_{n,E}^{-1}(\pi(q_0))$ and 
$C_r$ is any irreducible component of $C_2\times_{\P^1_2} E_r$.
Any point $q_r\in E_r$ such that $q_r\notin \Ram(\pi_r)$ 
(that is, its image in $\P^1_2$ is distinct from $r_1,...,r_4$) 
has the claimed property.   
\end{proof}

\begin{rem}
\label{rem:any}
Iterating this procedure (and adding isogenies) 
we obtain many elliptic curves $E'$ which are dominated by 
curves having an unramified cover onto $E$. 
It would be interesting to know if for any  
two elliptic curves over $\ovl{\Q}$ there exists 
a cycle connecting them (at least modulo isogenies). 
We will now show that {\em  any} elliptic curve can be
connected in this way to $E_0$.   
\end{rem}

Let  $E_0\subset \P^2=\{ (x:y:z)\}$ be the elliptic curve
$$
x^3+y^3+z^3=0,
$$
and 
$$
E_0[3]=\sT:=\left\{ 
\begin{array}{ccc} 
(1:0:1), & (1:0:-\zeta), & (1:0:-\zeta^2),\\           
(0:1:1), & (0:1:-\zeta), & (0:1:-\zeta^2),\\ 
(1:1:0), & (1:-\zeta:0), & (1:-\zeta^2:0)
\end{array} \right\}
$$ 
its set of $3$-torsion points (where $\zeta$ is a primitive cubic root of 1).
Denote by $\cE_{\la}=\{E_{\la}\}$ 
the family of elliptic curves on $\P^2$ passing through $\sT$ given by
$$
E_{\la}\,\,:\,\, x^3 + y^3 +z^3 + \la xyz = 0.
$$
It is easy to see that for each $\la$ the set $E_\la[3]$ of 3-torsion
points of $E_{\la}$ is precisely $\sT$.   
$$
\begin{array}{ccccc}
\pi & : & \P^2 & \ra &  \P^1\\
    &   &  (x:y:z) & \mapsto & (x+z:y)
\end{array}  
$$
the projection respecting the involution $x\ra z$ on $\P^2$.
Denote by $\pi_{\la}$ the restriction of $\pi$ to $E_{\la}$.
Clearly, $\pi_{\la}$ exhibits each $E_{\la}$ as a double cover of $\P^1$
and $\pi_{\la}$ has only simple double points for all $\la$. 
Moreover,
$$
\pi(\sT)=\{ (0:1),\, (1:-\zeta),\, (1:- \zeta^2),\, (1:-1),\, (1: 0)\}
$$
and for all $\la$ there exists a (non-empty)
set $S_{\la}\subset \Bran(\pi_{\la})\subset \P^1$ 
such that $\pi_{\la}^{-1}(S_{\la})\subset \sT$. 
Let $\pi_0'\,:\, E_0'\ra \P^1$ 
be a double cover ramified in 4 points in $\pi(\sT)$.

\

\begin{lemm}
\label{lemm:e}
Let $\iota\,:\, C\ra E_{\la}$ 
be a double cover such that over at least one point in $\Bran(\iota)$ 
the local ramification indices are even. 
Then there exists an unramified cover $\tilde{C}\ra C$ and a 
surjective morphism $\tilde{\iota}\,:\, \tilde{C}\ra E_0'$
such that over at least one point in $\Bran(\tilde{\iota})\subset E_0'$ 
all local ramification indices of $\tilde{\iota}$ are even. 
\end{lemm}

\begin{proof}
Consider the diagram 

\

\centerline{
\xymatrix{
E_{\la}\ar[d]_{\varphi_3}  & \ar[l]^{\iota}\ar[d]C_1 \\
E_{\la} \ar[d]_{\pi_{\la}} & \ar[l]  C \\
\P^1           &   
   }}

\

Then $C_1\ra \P^1$ has even local ramification 
indices over all points in $\pi(\sT)$. 
It follows that
$$
\tilde{C}:=C_1\times_{\P^1}E_0'\ra E_0'
$$
has even local ramification indices over the preimages of the fifth point
in $\pi(\sT)$, as claimed. 
\end{proof}

\begin{nota}
\label{nota:uni}
Let $\cC$ be the class of curves such that
there exists an elliptic curve $E$, a surjective map
$\iota\,:\, C\ra E$ and a point $q\in \Bran(\iota)$ such that 
all local ramification indices at points in $\iota^{-1}(q)$ are even.
\end{nota}

\begin{exam}
Any hyperelliptic curve of genus $\ge 2 $ belongs to $\cC$. More generally, 
$\cC$ contains any curve $C$ admitting a map  
$C\ra \P^1$ with even local ramification indices over
at least 5 points in $\P^1$. 
\end{exam}

\begin{prop}
\label{prop:ccc}
For any $C\in \cC$ there exists an unramified cover $\tilde{C}\ra C$ 
surjecting onto $C_0$ (with $C_0\ra \P^1$ given by 
$\sqrt[6]{z(1-z)}$). 
\end{prop}

\begin{proof}
Consider $C_1=C\in \cC$ with $\iota_1\,:\, C_1\ra E=E_{\la}$ 
as in \ref{nota:uni}. Define $C_2$ as an irreducible 
component of  $C_1\times_E E$:

\

\centerline{
\xymatrix{
C_1\ar[d]_{\iota_1} & \ar[l]_{\tau_2}\ar[d]^{\iota_2} C_2 \\
E               &  \ar[l]^{\varphi_3} \ar[d]^{\pi_\la}  E  \\ 
                                   &  \P^1    }}

\

\noindent
Define $C_3:=C_2\times_{\P^1} E_0$ by  the diagram

\

\centerline{
\xymatrix{
C_2\ar[d]_{\sigma_2} & \ar[l]_{\tau_3}\ar[d]^{\iota_3}C_3 \\
\P^1             & \ar[l]^{\pi_0} E_{0}.  
   }}

\

\noindent
Observe that for $q\in \Bran(\pi_0)$ the local ramification indices
in the preimage $(\iota_2\circ \pi_{\la})^{-1}(q)$ are all even. 
It follows that the map 
$\tau_3\,:\, C_3\ra C_2$ is {\em unramified} and 
that $\iota_3\, :\, C_3\ra E_0$ has even local ramification indices
over (the preimage of) 
$q_5\in \{\pi(\sT)\setminus \Bran(\pi_0)\}$ (the 5th point). 
Define $C_4$ as an irreducible component of 
$C_3\times_{E_0} E_{0}$ in the diagram

\

\centerline{
\xymatrix{
C_3\ar[d]_{\iota_3} & \ar[l]_{\tau_4}\ar[d]^{\iota_4}C_4 \\
E_0              & \ar[l]^{\varphi_{3}}E_{0}.  
   }}

\

\noindent
The map $\iota_4$ is ramified over the preimages 
$(\pi_0\circ \varphi_3)^{-1}(q_5)$, 
with even local ramification indices. 
Finally, $C_5= C_4\times_{E_0}  C_0$ from the diagram

\

\centerline{
\xymatrix{
C_4\ar[d]_{\iota_4}        &\ar[l]_{\tau_5}\ar[d] C_5 \\
E_0            & \ar[l]_{\iota_0}  C_{0}. 
   }}

\

\noindent
has a dominant map  onto $C_0$ and is unramified over $C_4$
(and consequently, $C_1$). 
\end{proof}

\end{document}